\documentclass[11pt]{article}
\usepackage{amssymb,amsmath,amsthm}

\newcommand{\pr}{\mbox{Pr}}
\newcommand{\con}{\mbox{Con}}

\date{}

\begin{document}

\title{The Surprise Examination Paradox and
the Second Incompleteness Theorem\thanks{First published in 
Notices of the AMS volume 57 number 11 (December 2010), published by
the American Mathematical Society.}}

\author{
Ran Raz\thanks{{\tt ran.raz@weizmann.ac.il}, Research supported by
the Israel Science Foundation (ISF),
the Binational Science Foundation (BSF)
and the Minerva Foundation.}\\
Weizmann Institute
}

\author{
Shira Kritchman \thanks{Faculty of Mathematics and
Computer Science, Weizmann Institute, Rehovot, Israel.
Email:~\texttt{shirrra@gmail.com}} \\
Weizmann Institute
\and Ran Raz  \thanks{Faculty of Mathematics and
Computer Science, Weizmann Institute, Rehovot, Israel.
Email:~\texttt{ran.raz@weizmann.ac.il}} \\
Weizmann Institute}

\maketitle

\begin{abstract}
We give a  new proof for G$\ddot{\mbox{o}}$del's second
incompleteness theorem, based on Kolmogorov complexity, Chaitin's
incompleteness theorem, and an argument that resembles the surprise examination paradox.

 We then go the other way around and suggest that the second incompleteness theorem gives a possible resolution of the surprise examination paradox. Roughly speaking, we argue that the flaw in the
derivation of the paradox is that it contains a hidden assumption that one can prove the consistency
of the mathematical theory in which the derivation is done; which is impossible by the second incompleteness theorem.
\end{abstract}

Few theorems in the history of mathematics have inspired
mathematicians and philosophers as much as G$\ddot{\mbox{o}}$del's
incompleteness theorems. The first incompleteness theorem states
that for any rich enough\footnote{We require that the theory can
express and prove basic arithmetical truths. In particular, ZFC
and Peano Arithmetic (PA) are rich enough.} consistent mathematical
theory\footnote{Here and below, we only consider first order theories with
recursively enumerable sets of axioms. For simplicity, let us
assume that the set of axioms is computable.}, there exists a
statement that cannot be proved or disproved within the theory.
The second incompleteness
theorem states that for any rich enough consistent mathematical theory, the
consistency of the theory itself cannot be proved (or disproved)
within the theory.

\section*{The First Incompleteness Theorem}

G$\ddot{\mbox{o}}$del's original proof for the first
incompleteness theorem~\cite{Godel} is based on the {\em liar paradox}.
\begin{quote}
{\bf The liar paradox:} consider the statement ``this statement is
false''. The statement can be neither true nor false.
\end{quote}
G$\ddot{\mbox{o}}$del considered the related statement ``this
statement has no proof''.
He showed that this statement
can be expressed in any theory that is capable of expressing
elementary arithmetic. If the statement has a proof, then it is
false; but since in a consistent theory any statement that has a
proof must be true, we conclude that if the theory is
consistent the statement has no proof. Since the
statement has no proof, it is true (over ${\mathbb N}$).
Thus, if the theory is consistent,
we have an example for a true statement
(over ${\mathbb N}$)
that has no proof.

The main conceptual difficulty in G$\ddot{\mbox{o}}$del's original
proof is the self-reference of the statement ``this statement
has no proof''. A conceptually simpler proof of the first
incompleteness theorem, based on  {\em Berry's paradox}, was given
by Chaitin~\cite{Chaitin}.
\begin{quote}
{\bf Berry's paradox:} consider the expression ``the smallest
positive integer not definable in under eleven words''. This
expression defines that integer in under eleven words.
\end{quote}

To formalize Berry's paradox, Chaitin uses the notion of {\em
Kolmogorov complexity}. The Kolmogorov complexity $K(x)$ of an
integer $x$ is defined to be the length (in bits) of the shortest
computer program that outputs $x$ (and stops). Formally, to define
$K(x)$ one has to fix a programming language, such as {\em LISP},
{\em Pascal} or {\em C++}. Alternatively, one can define $K(x)$ by
considering any {\em universal Turing machine}.

Chaitin's incompleteness theorem states that for any rich enough consistent
mathematical theory, there exists a (large enough) integer $L$
(depending on the theory and on the programming language that
is used to define Kolmogorov complexity), such that, for any
integer $x$, the statement ``$K(x)> L$'' cannot be proved within
the theory.

The proof given by Chaitin is as follows. Let $L$ be a large
enough integer. Assume for a contradiction that for some integer
$x$, there is a proof for the statement ``$K(x) > L$''. Let $w$ be
the first proof (say, according to the lexicographic order) for a
statement of the form ``$K(x) > L$''. Let $z$ be the integer $x$
such that $w$ proves ``$K(x) > L$''. It is easy to give a computer
program that outputs $z$: the program enumerates all possible
proofs $w$, one by one, and for the first $w$ that proves a
statement of the form ``$K(x) > L$'', the program outputs $x$ and
stops. The length of this program is a constant + $\log L$. Thus,
if $L$ is large enough, the Kolmogorov complexity of $z$ is less
than $L$. Since $w$ is a proof for ``$K(z)
> L$'' (which is a false statement), we conclude that the theory is
inconsistent.

Note that the number of computer programs of length $L$ bits is at
most $2^{L+1}$. Hence, for any integer $L$, there exists an integer $0
\leq x \leq 2^{L+1}$, such that $K(x)
> L$. Thus, for some integer $x$, the statement ``$K(x)
> L$'' is a true statement (over ${\mathbb N}$) that has no proof.

A different proof for G$\ddot{\mbox{o}}$del's first incompleteness
theorem,  also based on Berry's paradox, was given by
Boolos~\cite{Boolos} (see also~\cite{Vopenka,Kikuchi_Boolos}).
Other proofs for the first incompleteness theorem are also known (for a recent survey, see~\cite{Kotlarski_survey}).

\section*{The Second Incompleteness Theorem}

The  second incompleteness theorem follows directly from
G$\ddot{\mbox{o}}$del's original proof for the first
incompleteness theorem. As described above, G$\ddot{\mbox{o}}$del
expressed the statement ``this statement has no proof''
and showed that, if the theory is consistent, this is a true
statement (over ${\mathbb N}$)
that has no proof.
Informally, since the proof that this is a true
statement can be obtained within any rich enough theory,
such as Peano Arithmetic (PA) or ZFC,
if the consistency of the theory itself can
also be proved within the theory, then the statement can be proved
within the theory, which is a contradiction. Hence, if the theory
is rich enough, the consistency of the theory cannot be proved
within the theory.

Thus, the  second incompleteness theorem follows directly from
G$\ddot{\mbox{o}}$del's original proof for the first
incompleteness theorem. However, the second incompleteness theorem
doesn't follow from Chaitin's and Boolos'  simpler proofs for the
first incompleteness theorem. The problem is that these
proofs only show the existence of a true statement (over ${\mathbb N}$)
that has no proof, without giving an explicit example of such a
statement.

A different proof for the second incompleteness theorem, based on
Berry's paradox, was given by Kikuchi~\cite{Kikuchi_second}. This proof is model
theoretic, and seems to us somewhat less intuitive for people who
are less familiar with model theory.
For previous model theoretic proofs for the second incompleteness theorem
see~\cite{Kreisel} (see also~\cite{Smorinski}).

\section*{Our Approach}

We give a new proof for the second incompleteness theorem, based
on Chaitin's incompleteness theorem and an argument that resembles the {\em surprise
examination paradox}, (also known as the {\em unexpected hanging
paradox}).

\begin{quote}
{\bf The surprise examination paradox:} the teacher announces in
class: ``next week you are going to have an exam, but you will not
be able to know on which day of the week the exam is held until
that day''. The exam cannot be held on Friday, because otherwise,
the night before the students will know that the exam is going to
be held the next day. Hence, in the same way, the exam cannot be
held on Thursday. In the same way, the exam cannot be held on any
of the days of the week.
\end{quote}

Let $T$ be a (rich enough) mathematical theory, such as PA or
ZFC. For simplicity, the reader can assume that $T$ is ZFC, the
theory of all mathematics; thus, any mathematical proof, and in
particular any proof in this paper, is obtained within $T$.

Let $L$ be the integer guaranteed by Chaitin's incompleteness
theorem. Thus, for any integer~$x$, the statement ``$K(x)> L$''
cannot be proved (in the theory $T$), unless the theory is
inconsistent. Note, however, that for any integer $x$, such that,
$K(x) \leq  L$, there is a proof (in $T$) for the statement
``$K(x) \leq L$'', simply by giving the computer program of length
at most $L$ that outputs $x$ and stops, and by describing the
running of that computer program until it stops.


Let $m$ be the number of integers $0 \leq x \leq 2^{L+1}$, such
that, $K(x)> L$. (The number $m$ is analogous to the day of
the week on which the exam is held in the surprise examination
paradox). Recall that since the number of computer programs of
length $L$ bits is at most $2^{L+1}$, there exists at least one
integer $0 \leq x \leq 2^{L+1}$, such that, $K(x) > L$. Hence, $m \geq
1$.


Assume that $m=1$. Thus, there exists a single integer $x \in \{0,
\ldots ,2^{L+1}\}$ such that $K(x)> L$, and every other integer $y \in
\{0, \ldots ,2^{L+1}\}$ satisfies $K(y) \leq  L$. In this case, one
can prove that $x$ satisfies $K(x)> L$ by proving that every other
integer $y \in \{0, \ldots ,2^{L+1}\}$ satisfies $K(y) \leq  L$ (and
recall that there is a proof for every such statement). Since we
proved that $m \geq 1$, the only $x$ for which we didn't
prove $K(x) \leq  L$ must satisfy $K(x)> L$.

Thus, if $m=1$ then for some integer $x$, the statement ``$K(x)>
L$'' can be proved (in $T$). But we know that for any integer $x$,
the statement ``$K(x)> L$'' cannot be proved (in $T$), unless the
theory is inconsistent. Hence, if the theory is consistent, $m \geq 2$.
Since we assume that $T$ is a rich enough theory, we
can prove the last conclusion in $T$. That is, we can prove in $T$
that: if $T$ is consistent then $m \geq 2$.

Assume for a contradiction that the consistency of $T$ can be
proved within $T$. Thus,  we can prove in $T$ the statement ``$m
\geq 2$''. In the same way, we can work our way up and prove that
$m
\geq i+1$, for every $i \leq  2^{L+1} + 1$. In particular,
$m > 2^{L+1} + 1$,
which is a contradiction, since $m \leq 2^{L+1} + 1$
(by the definition of $m$).

\section*{The Formal Proof}

To present the proof formally, one needs to be able to
express provability within $T$, in the language of $T$. The
standard way of doing that is by assuming that the language of $T$
contains the language of arithmetics and by encoding every formula and
every proof in $T$ by an integer, usually referred to as the
G$\ddot{\mbox{o}}$del number of that formula or proof. For a
formula $A$, let $\ulcorner A \urcorner$ be its
G$\ddot{\mbox{o}}$del number. Let $\pr_T(\ulcorner A \urcorner)$
be the following formula: {\em there exists $w$ that is the
G$\ddot{\mbox{o}}$del number of a $T$-proof for the formula $A$}.
Intuitively, $\pr_T(\ulcorner A \urcorner)$ expresses the
provability of the  formula $A$. Formally, the formulas
$\pr_T(\ulcorner A \urcorner)$ satisfy the so-called
Hilbert-Bernays derivability conditions (see, for example,~\cite{Mendelson}):

\begin{enumerate}

\item If $T$ proves $A$ then $T$ proves $\pr_T(\ulcorner A \urcorner)$.

\item $T$ proves:
$\pr_T(\ulcorner A \urcorner) \rightarrow
\pr_T(\ulcorner\pr_T(\ulcorner A \urcorner)\urcorner)$.

\item $T$ proves:
$\pr_T(\ulcorner A \rightarrow B \urcorner) \rightarrow
\left(\pr_T(\ulcorner A \urcorner)
\rightarrow \pr_T(\ulcorner B \urcorner)\right)$

\end{enumerate}
The consistency of $T$ is usually expressed as the formula
$\con(T) \equiv \neg \pr_T(\ulcorner 0=1 \urcorner)$. In all that
comes below, $T \vdash A$ denotes ``$T$ proves $A$''.
We will prove that $ T \not \vdash \con(T)$, unless $T$ is inconsistent.

For our  proof, we will need two facts about provability
of claims concerning Kolmogorov complexity. First, we
need to know that $\con(T) \rightarrow
\neg \pr_T(\ulcorner K(x) > L \urcorner)$. We will use the following form
of Chaitin's incompleteness theorem
(see, for example,~\cite{Kikuchi_second}, Theorem 3.3).
\begin{eqnarray} \label{e1}
T \vdash \con(T) \rightarrow \forall x \in \{0, \ldots ,2^{L+1}\} \:
\neg \pr_T(\ulcorner K(x) > L \urcorner)
\end{eqnarray}
Second, we need to know that $(K(y) \leq L) \rightarrow
\pr_T(\ulcorner K(y) \leq L \urcorner)$. We will use the following
form (formally, this follows since $K(y) \leq L$ is a $\Sigma_1$ formula;
see, for example,~\cite{Kikuchi_second}, Theorem~1.2 and Section~2).
\begin{eqnarray} \label{e2}
T \vdash  \forall y \in \{0, \ldots ,2^{L+1}\} \:
\left((K(y) \leq L) \rightarrow
\pr_T(\ulcorner K(y) \leq L \urcorner)\right)
\end{eqnarray}

Assume for a contradiction that $T$ is consistent and
$ T \vdash \con(T)$.
Then, by Equation~\ref{e1},
\begin{eqnarray} \label{e3}
 T \vdash \forall x \in \{0, \ldots ,2^{L+1}\} \;  \neg
\pr_T(\ulcorner K(x) > L \urcorner)
\end{eqnarray}

We will derive a contradiction by proving by induction
that, for every $i \leq 2^{L+1}+1$,
$ T \vdash (m \geq i+1)$, where $m$ is defined as in the previous section.
Since obviously $ T \vdash (m \leq 2^{L+1}+1)$,
this is a contradiction to the assumption
that $T$ is consistent and $ T \vdash \con(T)$.
Since we already know that $ T \vdash (m \geq 1)$, we already have the
base case of the induction. Assume
(the induction hypothesis)
that for some $1 \leq i \leq 2^{L+1}+1$,
\begin{eqnarray*}
T \vdash (m \geq i )
\end{eqnarray*}
We will show that $ T \vdash (m \geq i+1)$, as follows.
Let $r=2^{L+1}+1-i$.

\begin{enumerate}
\item
By the definition of $m$,

$ T \vdash (m =i)
\rightarrow
\exists \; \mbox{different} \;
y_1, \ldots, y_r \in \{0, \ldots ,2^{L+1}\}\;
\bigwedge_{j=1}^r (K(y_j) \leq L )$

\item
Hence, by Equation~\ref{e2},

$ T \vdash (m =i)
\rightarrow
\exists \; \mbox{different} \;  y_1, \ldots, y_r \in
\{0, \ldots ,2^{L+1}\} \;
\bigwedge_{j=1}^r \pr_T(\ulcorner K(y_j) \leq L \urcorner)$

\item
For every
$\mbox{different} \;  y_1, \ldots, y_r \in \{0, \ldots ,2^{L+1}\}$,
and every $x \in \{0, \ldots ,2^{L+1}\}
\setminus \{ y_1,\ldots,y_{r}  \}$,

$ T \vdash ( m \geq i ) \rightarrow
\left(\bigwedge_{j=1}^r (K(y_j) \leq L )
\rightarrow (K(x) > L)\right),
 $

(by the definition of $m$), and
hence by Hilbert-Bernays derivability conditions,

$ T \vdash \pr_T(\ulcorner m \geq i \urcorner) \rightarrow
\left(\bigwedge_{j=1}^r \pr_T(\ulcorner K(y_j) \leq L \urcorner)
\rightarrow \pr_T(\ulcorner K(x) > L \urcorner)\right)
 $

\item
By the previous two items,

$ T \vdash ((m =i) \wedge \pr_T(\ulcorner m \geq i \urcorner))
\rightarrow
\exists x \in \{0, \ldots ,2^{L+1}\} \; \pr_T(\ulcorner
K(x) > L \urcorner) $

\item
Since $ T \vdash (m \geq i )$ (by the induction hypothesis),
$ T \vdash \pr_T(\ulcorner m \geq i \urcorner)$.
Hence,

$ T \vdash (m =i)
\rightarrow
\exists x \in \{0, \ldots ,2^{L+1}\} \;
\pr_T(\ulcorner K(x) > L \urcorner) $

\item
Hence, by Equation~\ref{e3},

$ T \vdash \neg (m =i)$

\item
Hence, since $ T \vdash (m \geq i )$ ,

$ T \vdash (m \geq i+1)$
\end{enumerate}
\qed

\section*{A Possible Resolution of The Surprise Examination Paradox}

In the previous sections we gave a proof for G$\ddot{\mbox{o}}$del's second incompleteness theorem by an argument that resembles the surprise examination paradox.
In this section we go the other way around and suggest that the second incompleteness theorem gives a possible resolution of the surprise examination paradox. Roughly speaking, we argue that the flaw in the
derivation of the paradox is that it contains a hidden assumption that one can prove the consistency
of the mathematical theory in which the derivation is done; which is impossible by the second incompleteness theorem.

The important step in analyzing the paradox is the translation of the
teacher's announcement into a mathematical language. The key point lies
in the formalization of the notions of surprise and knowledge.

As before, let $T$
be a rich enough mathematical theory (say, ZFC). Let $\left\{ 1,\dots,5\right\} $ be the days of the week and let $m$ denote the day of the week on which the
exam is held.
Recall the teacher's
announcement: {}``next week you are going to have an exam, but you
will not be able to know on which day of the week the exam is held
until that day''.  The first part of the announcement is formalized as $m\in\left\{ 1,\dots,5\right\} $.
A standard way that appears in the literature to formalize the second part is
by replacing the notion of {\bf knowledge} by the notion of {\bf provability}~\cite{Shaw,Fitch} (for a recent survey see~\cite{Chow}). The second part is rephrased as {}``on the night before
the exam you will not be able to prove, using this statement, that
the exam is tomorrow'', or, equivalently,
``for every $1 \le i \le 5$, if you are able to prove, using this statement, that $(m \ge i) \rightarrow (m=i)$, then $m \ne i$''.
This can be formalized as the following statement that we denote by $S$ (the statement $S$ contains both parts of the teacher's announcement):

\begin{equation*} \label{e4}
S \;\; \equiv \;\; [m\in\left\{ 1,\dots,5\right\}] \bigwedge_{1\le i\le5}\left[\pr_{T,S}\left(\ulcorner m\ge i\to m=i\urcorner\right)\to (m\neq i)\right]
\end{equation*}
where $\pr_{T,S}(\ulcorner A\urcorner)$ expresses the provability of a formula $A$ from the formula $S$ in the theory $T$,
(formally, $\pr_{T,S}(\ulcorner A\urcorner)$ is the formula: {\em there exists $w$ that is the
G$\ddot{\mbox{o}}$del number of a $T$-proof for the formula $A$ from the formula $S$}).
Note that the formula $S$ is self-referential. Nevertheless, it is well known that this is not a real problem and that such a formula $S$ can be formulated (see~\cite{Shaw,Chow}; for more about this issue, see below).

Let us try to analyze the paradox when the teacher's announcement is formalized as the above statement $S$.
We will start from the last day.
The statement
$m\ge 5$ together with $m\in\left\{ 1,\dots,5\right\} $ imply $m=5$. Hence,
$\pr_{T,S}\left(\ulcorner m\ge 5\to m=5\urcorner\right)$ and by $S$ we can
conclude $m \ne 5$. Thus, $S$ implies   $m\in\left\{ 1,\dots,4\right\} $. In the same way, working our way down, we can prove
$\pr_{T,S}\left(\ulcorner m\ge 4\to m=4\urcorner\right)$ and by $S$
 we can conclude   $m \ne 4$. In the same way, $m \ne 3$, $m \ne 2$, and $m \ne 1$.
In other words,  $S$ implies
$m\notin\left\{ 1,\dots,5\right\} $.
Thus, $S$ contradicts itself.

The fact that $S$ contradicts itself gives a certain explanation for the paradox; the teacher's announcement is just a contradiction.
On the other hand, we feel that this formulation doesn't fully explain the paradox: Note that since $S$ is a contradiction it can be used to prove any statement. So, for example, on Tuesday night the students can use $S$ to  prove that the exam will be held on Wednesday. Is it fair to say that this means that they {\bf know} that the exam will be held on Wednesday? No, because they can also use $S$ to prove that the exam will be held on Thursday. Thus, we conclude that since $S$ is a contradiction, {\bf provability} from $S$ doesn't imply {\bf knowledge}. Recall, however, that the very intuition behind the formalization of the teacher's announcement as $S$ was that the notion of {\it knowledge} can be replaced by the notion of {\it provability}. But if provability from $S$ doesn't imply knowledge, the statement $S$ doesn't seem to be an accurate translation of the teacher's announcement into a mathematical language.

Is there a better way to formalize the teacher's announcement? To answer this question, let us analyze the situation from the students' point of view on Tuesday night. There are three possibilities:

\begin{enumerate}
\item
On Tuesday night, the students are not able to prove that the exam will be held on Wednesday.
\item
On Tuesday night, the students are able to prove that the exam will be held on Wednesday, but they are also able to prove for some other day that the exam will be held on that day. \\
(Note that this possibility can only occur if the system is inconsistent, and is in fact equivalent to the inconsistency of the system).
\item
On Tuesday night, the students are able to prove that the exam will be held on Wednesday, and they are not able to prove for any other day that the exam will be held on that day.
\end{enumerate}
We feel that only in the third case is it fair to say that the students know that the exam will be held on Wednesday.
They know that the exam will be held on Wednesday only if
they are able to prove that the exam will be held on Wednesday, and they are not able to prove for any other day that the exam will be held on that day.

We hence rephrase the second part of the teacher's announcement as
``for every $1 \le i \le 5$, if one can prove (using this statement) that $(m \ge i) \rightarrow (m=i)$, and there is no $j \ne i$ for which one can prove (using this statement) $(m \ge i) \rightarrow (m =j)$,
then $m \ne i$''.
Thus, the teacher's announcement is the following statement\footnote{This statement is equivalent to one of the suggestions (the statement $I_5$) made by Halpern and Moses~\cite{HM}. However, the analysis of the paradox there is different from the one shown here and makes no use of G$\ddot{\mbox{o}}$del's second
incompleteness theorem.}:

\begin{equation*}
S \;\; \equiv \;\; [m\in\left\{ 1,\dots,5\right\}] \bigwedge_{1\le i\le5}\left[\left(\pr_{T,S}\left(\ulcorner m\ge i\to m=i\urcorner\right)\bigwedge_{1\le j\le5,j\neq i}\neg\pr_{T,S}\left(\ulcorner m\ge i\to m=j\urcorner\right)\right)\to (m\neq i)\right]
\end{equation*}

Let us try to analyze the paradox when the teacher's announcement is formalized as the new statement $S$.
As before, $m\ge 5$ together with $m\in\left\{ 1,\dots,5\right\} $ imply $m=5$. Hence,
$\pr_{T,S}\left(\ulcorner m\ge 5\to m=5\urcorner\right)$.
However, this time one cannot use  $S$ to
conclude $m \ne 5$, since it is possible that for some  $j \ne 5$ we also have
$\pr_{T,S}\left(\ulcorner m\ge 5\to m=j\urcorner\right)$.
This happens iff the system $T + S$ is inconsistent.
Formally, this time one cannot use $S$ to deduce $m \ne 5$, but rather the formula  $$\con(T,S) \rightarrow (m \ne 5),$$ where $\con(T,S) \equiv \neg \pr_{T,S}(\ulcorner 0=1 \urcorner)$ expresses
the consistency of $T+S$.  Since by the second incompleteness theorem one cannot prove $\con(T,S)$ within $T+S$, we cannot conclude that $S$ implies $m \ne 5$ and hence cannot continue the argument.

More precisely, since $S$ doesn't imply $m \in \{1,\ldots,4\}$, but rather $\con(T,S) \rightarrow m \in \{1,\ldots,4\},$ when we try to work our way down we do not get the desired formula
$\pr_{T,S}\left(\ulcorner m\ge 4\to m=4\urcorner\right)$, but rather the formula  $$\pr_{T,S}\left(\ulcorner \con(T,S) \wedge (m\ge 4) \to m=4\urcorner\right),$$
which is not enough to continue the argument.

Thus, our conclusion is that if the students believe in the consistency of $T+S$ the exam cannot be held on Friday, because on Thursday night the students will know that if $T+S$ is consistent the exam will be held on Friday. However, the exam can be held on any other day of the week because the students cannot prove the consistency of $T+S$.

Finally, for completeness, let us address the issue of the self-reference of the statement $S$. The issue of self-referentiality of a statement goes back to G$\ddot{\mbox{o}}$del's original proof for the first
incompleteness theorem. The self-reference is what makes G$\ddot{\mbox{o}}$del's original proof conceptually difficult, and what makes the teacher's announcement in the surprise examination paradox paradoxical.

To solve this issue, G$\ddot{\mbox{o}}$del introduced the
technique of diagonalization. The same
technique can be used here.
To formalize $S$, we will use the notation $a \Rightarrow b$ to indicate implication between G$\ddot{\mbox{o}}$del numbers $a$ and $b$. That is,
$a \Rightarrow b$
is a statement indicating that $a$ is a G$\ddot{\mbox{o}}$del
number of a statement $A$, and $b$ is a G$\ddot{\mbox{o}}$del
number of a statement $B$, such that, $A \rightarrow B$.
We will also need the function $\mbox{Sub}(a,b)$ that represents substitution of $b$ in the formula with G$\ddot{\mbox{o}}$del number~$a$. That is, if $a$ is a G$\ddot{\mbox{o}}$del number of a formula $A(x)$ with free variable $x$, and $b$ is a number, then $\mbox{Sub}(a,b)$ is the G$\ddot{\mbox{o}}$del number of the statement $A(b)$.

Let $v_{ij} \equiv \ulcorner m\ge i\to m=j\urcorner$.
Denote by $Q(x)$ the formula
$$
 Q(x)\;\; \equiv \;\; [m\in\left\{ 1,\dots,5\right\}] \bigwedge_{1\le i\le5}
 \left[\left(\pr_{T}\left(\mbox{Sub}(x,x) \Rightarrow v_{ii}\right)\bigwedge_{1\le j\le5,j\neq i}\neg\pr_{T}\left(\mbox{Sub}(x,x) \Rightarrow v_{ij}\right)\right)\to (m\neq i)\right]
$$
Let $q$ be the G$\ddot{\mbox{o}}$del number of the formula $Q(x)$. The statement $S$ is formalized as  $S \equiv Q(q)$.
To see that this statement is the one that we are interested in, denote by $s$ the G$\ddot{\mbox{o}}$del number of $S$ and note that $s = \mbox{Sub}(q,q)$. Thus,
$$
 S \;\; \equiv \;\; [m\in\left\{ 1,\dots,5\right\}] \bigwedge_{1\le i\le5}
 \left[\left(\pr_{T}\left(s \Rightarrow v_{ii}\right)\bigwedge_{1\le j\le5,j\neq i}\neg\pr_{T}\left(s \Rightarrow v_{ij}\right)\right)\to (m\neq i)\right]
$$

\end{document}